\documentclass[12pt,reqno,papersize]{amsart}
\usepackage{amsmath,amssymb}
\usepackage[dvips]{graphicx,color,xcolor}
\usepackage{latexsym,amssymb}
\usepackage{mathrsfs}

\newtheorem{thm}{Theorem}[section]

\newtheorem{lem}[thm]{Lemma}
\newtheorem{rem}[thm]{Remark}
\newtheorem{notation}[thm]{Notation}

\numberwithin{equation}{section}
\numberwithin{thm}{section}

\begin{document}

\begin{center}\large \bf 
Self-similar solutions for compressible Navier-Stokes 
equations 
\end{center}

\footnote[0]
{
{\it Mathematics Subject Classification} (2010): 35C06, 76N10

{\it 
Keywords}: Compressible Navier-Stokes equations; self-similar solutions; vacuum. 

{\it Acknowledgements}: P. Germain was partially supported by the NSF grant DMS-1501019. 
T. Iwabuchi was supported 
by JSPS Program for Advancing Strategic International Networks to Accelerate 
the Circulation of Talented Researchers 
and by JSPS KAKENHI Grant 17H04824. 

{\it Addresses}: Pierre Germain, Courant Institute of Mathematical Sciences, New York University, 251 Mercer Street, New York, N.Y. 10012-1185, USA, {\tt pgermain@cims.nyu.edu}

Tsukasa Iwabuchi, Mathematical Institute, Tohoku University, 
Sendai, 980-8578, JAPAN, 
{\tt t-iwabuchi@m.tohoku.ac.jp}
}

\vskip5mm

\begin{center}

Pierre GERMAIN,\quad Tsukasa IWABUCHI

\end{center}

\vskip5mm

\begin{center}
\begin{minipage}{135mm}
\footnotesize
{\sc Abstract. } We construct forward self-similar solutions (expanders) for the compressible Navier-Stokes equations. Some of these self-similar solutions are smooth, while others exhibit a singularity do to cavitation at the origin.

\end{minipage}
\end{center}

\tableofcontents


\section{Introduction}

\subsection{The model}

We study self-similar solutions of 
the following compressible Navier-Stokes equations 
in $\mathbb R^d$ with $d \geq 1$. 
\begin{equation}\label{eq:cNS_full}
\begin{cases}
\partial_t \rho + {\rm div \, } \big( \rho u \big) 
= 0 , 
& \quad t > 0, x \in \mathbb R^d, 
\\
\displaystyle 
\partial _t (\rho u) + 
{\rm div \,} \big( \rho u \otimes u \big) + \nabla \pi  
 = {\rm div } \, \tau,
& \quad t > 0, x \in \mathbb R^d, 
\\
\displaystyle 
\partial_t 
\Big[ \rho \Big( \frac{|u|^2}{2} + e \Big) \Big] 
+ {\rm div} 
\Big[ u \Big( \rho \Big(\frac{|u|^2}{2} + e \Big) 
      + \pi \Big) \Big]
- {\rm div } \, q = {\rm div \,} (\tau \cdot u) , 
& \quad t > 0, x \in \mathbb R^d, 
\end{cases}
\end{equation}
where $\rho(t,x)$ is the density of the fluid, $u(t,x)$ its velocity, $e(t,x)$ its internal energy, $\pi(t,x)$ its pressure, $\tau(t,x)$ its stress tensor, and finally $q(t,x)$ its internal energy flux. The fluid will furthermore be described by its temperature $\theta(t,x)$. 

We assume the following constitutive relations:
\begin{itemize}
\item Joule's first law: $$e = C_V \theta,$$ where $C_V > 0$ is the heat constant.
\item Ideal gas law: $$\pi = \rho R \theta,$$ where $R>0$ is the ideal gas constant.
\item Newtonian fluid: this implies 
$$\tau := \lambda {\rm div \,} u \, {\rm Id} 
 + 2 \mu D(u), 
\quad 
 D(u) = \frac{\nabla u + (\nabla u)^T}{2} , 
\quad \nabla u = (\partial_{x_i} u_j), 
$$
where $\lambda$ and $\mu$ are the Lam\'e coefficients, which satisfy
$$
\mu>0 \qquad \mbox{and} \qquad 2\mu + d\lambda \geq 0.
$$
\item Fourier's law $$q = \kappa \nabla \theta,$$ where $\kappa>0$ is the thermal conductivity.
\end{itemize}

We refer to~\cite{Lions} for a more detailed discussion of these assumptions. The equations become

\begin{equation}\label{eq:cNS} \tag{cNS}
\begin{cases}
\partial_t \rho + {\rm div \, } \big( \rho u \big) 
= 0 , 
\\
\displaystyle 
\partial _t (\rho u) + 
{\rm div \,} \big( \rho u \otimes u \big) + \nabla ( \rho R \theta ) 
 = (\lambda + \mu) \nabla {\rm div \,} u + \mu \Delta u,
\\
\displaystyle 
\partial_t 
\Big[ \rho \Big( \frac{|u|^2}{2} + C_V \theta \Big) \Big] 
+ {\rm div} 
\Big[ u \Big( \rho \Big(\frac{|u|^2}{2} + C_V \theta \Big) 
      + \rho R \theta \Big) \Big]
- \kappa \Delta \theta
\\
\qquad \qquad \qquad \qquad \qquad \qquad \qquad \qquad \qquad= {\rm div \,} (\lambda ({\rm div \,} u) u + 2 \mu D(u) \cdot u ) , 
\end{cases}
\end{equation}

\subsection{Forward self-similar solutions} The equations~\eqref{eq:cNS} exhibit a scaling invariance: the set of solutions is left invariant by the transformation
$$
\rho (t,x) \to \rho (\lambda ^2 t , \lambda x), \quad 
u (t,x) \to \lambda u(\lambda ^2 t , \lambda x), \quad 
\theta (t,x) \to \lambda^2 \theta(\lambda^2 t, \lambda x), \quad \mbox{for $\lambda>0$}.
$$
This scaling invariance suggests looking for self-similar solutions of the form
\begin{equation}\label{eq:selfsimilar}
\rho (t,x) = P \left( \frac{r}{\sqrt t} \right),  
\quad 
u(t,x) = \frac{1}{\sqrt{t}} U\left( \frac{r}{\sqrt t} \right)
  \frac{x}{r},
\quad 
\theta (t,x) = \frac{1}{t} \Theta \left( \frac{r}{\sqrt t} \right),
\end{equation}
where $r = |x|$, $P$, $U$ and $\Theta$ are 
scalar functions from $(0,\infty)$ to $\mathbb R$. 

It is natural to expect that there exist real numbers $P_\infty$, $U_\infty$, $\Theta_\infty$ such that
$$
P(r) \to P_\infty, \quad U(r) \sim \frac{U_\infty}{r}, \quad \Theta(r) \sim \frac{\Theta_\infty}{r^2} \quad \mbox{as $r \to \infty$},
$$
in which case this self-similar solution is associated to self-similar data
\begin{equation}
\label{initialdata}
(\rho,u,\theta)(t=0) = \left( P_\infty, \, U_\infty \frac{x}{r^2}, \, \frac{\Theta_\infty}{r^2} \right).
\end{equation}

\subsection{Known results} 

\subsubsection{Weak and strong solutions of~\eqref{eq:cNS}} 
In the very rich existing literature, we mention weak, finite-energy  solutions by Lions~\cite{Lions-1998}, variational solutions by 
Feireisl-Novotn\'y-Petzeltov\'a~\cite{FNP-2001} (see also Feireisl~\cite{F-2004}), classical solutions with finite energy by 
Matsumura-Nishida~\cite{MN-1980} (see also
Huang-Li~\cite{HL-2018} with vacuum), solutions in Besov spaces with the interpolation index one by Danchin~\cite{D-2001} (see also Chikami-Danchin~\cite{CD-2015}).

However, the initial velocity and temperature in \eqref{initialdata} are homogeneous functions of degree 1 and 2, respectively, 
which therefore do not fit any of these frameworks. Indeed,
\[
u_0 \not \in L^2 (\mathbb R^d) \cup L^d (\mathbb R^d)
 \cup \dot B^{\frac{d}{p}-1}_{p,1} (\mathbb R^d), \quad 
\theta_0 \not \in L^1 (\mathbb R^d) \cup L^{\frac{d}{2}} (\mathbb R^d)
 \cup \dot B^{\frac{d}{p}-2}_{p,1} (\mathbb R^d).
\]

Let us now try and be more specific, and explain why the classical construction methods cannot apply. Regarding weak solutions, the obstruction is obviously that the data~\eqref{initialdata} has infinite energy. Regarding strong solutions, the main obstacle is that the linear problem does not lead to $\nabla u \in L^1 (0,T; L^\infty)$, which in turn prevents any control of the density in $L^\infty$. Even worse, it is actually the case that $\int_0^1 e^{s (2\mu + \lambda)\Delta} {\rm div \, } u_0 \, ds \not \in L^\infty (\mathbb R^d)$.

\subsubsection{Self-similar solutions of~\eqref{eq:cNS}}
There are only few results in this direction. 
Under a different scaling property from the parabolic type \eqref{eq:selfsimilar}, 
Qin-Su-Deng~\cite{QSD-2008} proved the non-existence of forward and backward self-similar solutions to the compressible Navier-Stokes equations in one dimension. Local energy of forward and backward self-similar solutions 
was also investigated in \cite{QSD-2008} 
but the total energy blows up 
at $t = 0$ and $t = T$, respectively, where $T$ is the given time 
appearing in  the definition of backward self-similar solutions. 
We also refer to related papers~\cite{GJ-2006} by Guo-Jiang 
(isothermal compressible Navier-Stokes equations) and 
Li-Chen-Xie~\cite{LCX-2013} (density-dependent viscosity).  

\subsubsection{The case of incompressible Navier-Stokes} This case is different in two respects. First, the ansatz which we chose above 
(radial velocity) 
is incompatible with incompressibility, in fact, the velocity is irrotational; 
therefore, it is not possible to reduce the problem to a one-dimensional one, as we shall do in the present article. Second, the existence of forward self-similar solutions is known since strong solutions can be built up from  small self-similar data: see for instance Cannone and Planchon~\cite{CP}, Chemin~\cite{Chemin} and Koch and Tataru~\cite{KT}. The case of large self-similar data was recently treated by Jia and Sverak~\cite{JS}, who could prove the existence of smooth self-similar solutions.

\subsubsection{Vacuum state} 
Few papers are known related to the vacuum. 
Xin~\cite{Xin-1998} found the blow-up solutions for the initial density with the compact support. 
Hoff and Smoller~\cite{HS-2001} considered 1D barotropic Navier-Stokes equations 
and showed non-formulation of vacuum state due to 
the persistency of the almost everywhere positivity of the density. 
Jang and Masmoudi~\cite{JM-2015} 
proved local in time well-posedness of the 3D compressible Euler equations under the barotropic condition with a physical vacuum. 
We also refer to \cite{JM-2015-2} for the overview about problems of vacuum state. 

\subsection{Obtained results}

We only sketch below our two main results, and refer to Theorems~\ref{theo:2} and \ref{theo:3} for complete statements.

\begin{thm}[Smooth self-similar solutions; simplified statement] 
\label{thm}
Let $d \geq 3$. 
There exists a family of smooth self-similar solutions of the form \eqref{eq:selfsimilar} correponding to data~\eqref{initialdata},
where $P_\infty>0$, $U_\infty <0$, and $\Theta_\infty>0$. 
The parameters $U_\infty$ and $\Theta_\infty$ have to be chosen sufficiently small, and the allowed values of $P_\infty, U_\infty, \Theta_\infty$ form a two-dimensional manifold.

The profiles $(P(r),U(r),\Theta(r))$ are smooth functions of $[0,\infty)$ such that
\begin{align*}
& \inf_r P(r) > 0, \qquad \sup_r P(r) <\infty,\\
& |U(r)| \lesssim \frac{r}{(1+r)^3}, \qquad |U'(r)| \lesssim  \frac{r}{(1+r)^3}, \\
& |\Theta(r)| \lesssim \frac{1}{(1+r)^2}, \qquad |\Theta' (r)| \lesssim \frac{r}{(1+r)^2}, 
\end{align*}
and furthermore
\begin{align*}
& P(r) = P_\infty + O \left( \frac{1}{r^2} \right), \\
& U(r) = \frac{U_\infty}{r} + O \left( \frac{1}{r^3} \right), \\
& \Theta(r) = \frac{\Theta_\infty}{r^2} + O\left( \frac{1}{r^4} \right). 
\end{align*}
\end{thm}

The previous theorem can be thought of as perturbative, around the trivial (self-similar) solution $(\rho,u,\theta) = (Constant,0,0)$.

\begin{thm}[Cavitating self-similar solutions; simplified statement]
Let $d \geq 3$. There exists a family of self-similar solutions of the form \eqref{eq:selfsimilar} correponding to data~\eqref{initialdata}, where $P_\infty > 0, U_\infty>0 , \Theta _ \infty > 0$. The parameter $P_\infty$ has to be chosen sufficiently small, and the allowed values of $P_\infty, U_\infty, \Theta_\infty$ form a three-dimensional set.

The profiles $(P(r),U(r),\Theta(r))$ are smooth functions of $(0,\infty)$, which, for $r \to 0$ behave as follows:
\begin{align*}
& P(r) = P_\delta \left( \frac{r}{\delta} \right)^{\frac{2d\alpha}{1-2\alpha}} + O \left( \frac{r}{\delta} \right)^{\frac{2d\alpha}{1-2\alpha}+1+d\alpha},\\
& U(r) = \alpha r + O(r^{1 + \frac{2d\alpha}{1-2\alpha}}), \\
& \Theta(r) = \Theta_0 + O(r^2),
\end{align*}
where $\alpha$ and $P_\delta$ are small parameters.

The profiles $(P(r),U(r),\Theta(r))$ also satisfy the global bounds
\begin{align*}
& |P(r)| \lesssim P_\delta \min\left[1 , \left( \frac{r}{\delta} \right)^{\frac{2d\alpha}{1-2\alpha}} \right], \\
& |U(r)| \lesssim \frac{\alpha r}{(1+ \sqrt{P_\delta} r)^2}, \qquad |U'(r)| \lesssim \frac{\alpha}{(1+\sqrt{P_\delta} r)^2}, \\
& |\Theta(r)| \lesssim \frac{1}{(1+\sqrt{P_\delta}r)^2}, \qquad |\Theta'(r)| \lesssim \frac{\sqrt{P_\delta} r}{(1+\sqrt{P_\delta}r)^2}.
\end{align*}
and finally
\begin{align*}
& P(r) = P_\infty + O \left( \frac{1}{r^2} \right), \\
& U(r) = \frac{U_\infty}{r} + O \left( \frac{1}{r^3} \right), \\
& \Theta(r) = \frac{\Theta_\infty}{r^2} + O\left( \frac{1}{r^4} \right).
\end{align*}
\end{thm}

\begin{rem} \begin{itemize}
\item If $d = 1,2$, solutions can be constructed in a very similar way to the above theorems. However, we excluded $d=1,2$ because an initial data of the type $\frac{\Theta_\infty}{r^2}$ is not locally integrable, and thus does not make sense in the sense of distributions. Furthermore, we could not ensure positivity of $\Theta$.
\item Although $u (t,x) := t^{-1/2} U(t^{-1/2}|x|)x / |x| \notin L^1 (0,1 ; Lip(\mathbb R^d))$, one can define Lagrangian coordinates for the velocity fields defined in the two above theorems.
\end{itemize}
\end{rem}
\vskip3mm 

\subsection{Organization of the paper}
In Section 2, we derive the integro-differential equations which result from our ansatz.

In Section 3, we state a complete version of the existence theorem in the smooth case, and proceed to prove it.

In Section 4, we state a complete version of the existence theorem in the cavitating case, and proceed to prove it.

\section{ODEs and integro-differential equations} 

\subsection{Derivation of the system of ODEs}

Consider solutions such that \eqref{eq:selfsimilar} is satisfied. 
Let us starting by proving that 
the partial differential equations \eqref{eq:cNS} is 
equivalent to the following ordinary differential equations for any 
$r = |x| > 0$: 
\begin{equation}\label{eq:ODEs}
\begin{cases}
\displaystyle 
-\frac{1}{2} r P ' + P' U 
+ P \Big( U' + \frac{d-1}{r} U \Big) 
= 0 , 
\\[5mm]
\begin{split}
\displaystyle 
-\frac{1}{2} 
& P  U 
- \frac{1}{2} r (P U)' 
+ (P U ^2) ' + \frac{d-1}{r} P U^2 
+(P R \Theta)'
\\
= 
& 
(2\mu + \lambda) 
 \Big( U'' + \frac{d-1}{r} U' - \frac{d-1}{r^2} U \Big) ,
\end{split}
\\[10mm]
\displaystyle 
\begin{split}
-P 
& \Big( \frac{U^2}{2}
  + C_V \Theta \Big) 
- \frac{1}{2} r 
  \Big( P \Big( \frac{U^2}{2} + C_V \Theta \Big) \Big) '
+ \Big( U P \Big( \frac{U^2}{2} + C_V \Theta \Big) 
        + U P R \Theta \Big) '
\\
&
+ \frac{d-1}{r} 
  \Big( U P \Big( \frac{U^2}{2} + C_V \Theta \Big)
       + U P R \Theta
  \Big)
- \kappa \Big( \Theta '' + \frac{d-1}{r} \Theta ' \Big)
\\
= 
& 
2\mu \Big( (U')^2 + \frac{d-1}{r^2} U^2 \Big) 
+ \lambda \Big( U' + \frac{d-1}{r} U \Big) ^2
\\
& 
+ (2\mu + \lambda) 
 \Big( U'' + \frac{d-1}{r}U' - \frac{d-1}{r^2} U \Big) U. 
\end{split}
\end{cases}
\end{equation}

\bigskip

The above equations follow in a straightforward manner from the formulas ($f$ denoting a scalar function)
\begin{align*}
& \nabla f(r) = f'(r) \frac{x}{r}, \\
& {\rm div} \left( f(r) \frac{x}{r} \right) = f' + \frac{d-1}{r} f, \\
& {\rm div} \left( f(r) \frac{x}{r} \otimes \frac{x}{r} \right)= \left( f' + \frac{d-1}{r} f \right) \frac{x}{r}, \\
& \Delta  \left( f(r) \frac{x}{r} \right) = \left( f'' + \frac{d-1}{r} f'(r) - \frac{d-1}{r^2} f(r) \right) \frac{x}{r}, \\
& D\left( f(r) \frac{x}{r} \right) = \frac{f(r)}{r} \operatorname{Id} + \left( f'(r) - \frac{f(r)}{r} \right) \frac{x}{r} \otimes \frac{x}{r}.
\end{align*}

\subsection{Integro-differential equations for smooth solutions}
We next write the ordinary differential equations \eqref{eq:ODEs} as integral equations under the condition at 
$r = 0$ that 
\begin{align*}
& P(0) = P_0>0, \qquad P'(r) = O(1), \\
& U(r) = O(r^2), \qquad U'(r) = O(r)  \qquad 
 \text{with } \frac{1}{2} r - U(r) >0 \text{ for any } r>0,
\\
& \Theta(0) = \Theta_0>0, \qquad \Theta'(r) =O(r).
\end{align*}
We will obtain the following formulas: 
\begin{align}\label{0728-1}
P (r) = & e^{V(r)} P_0, 
\\ \label{0728-2}
U(r) 
=
& \frac{r^{-d+1}}{2\mu + \lambda} 
\int_0^r 
  r_1 ^{d-1} e^{-W(r) + W(r_1)} F_U(r_1) \, dr_1,
\\ \label{0728-3}
\Theta (r) 
= 
&
(d-2)r^{-d+2}
\int_0^r r_1 ^{d-3} e^{-Z(r) + Z (r_1)} \, dr_1 \Theta_0
- \frac{U^2}{2C_V}
+ \frac{r^{-d+2}}{\kappa} \int_0^r 
  r_1^{d-2} e^{-Z(r) + Z(r_1)} F_\Theta (r_1)
 \, dr_1 ,
\end{align}
where 
\begin{equation} \notag 
\begin{split}
F_U(r_1)
:=
&      P U^2 
       + \int_0^{r_1} \frac{d-1}{r_2} P U^2 \, dr_2 
       + P R \Theta - P_0 R \Theta_0, 
\\
F_\Theta (r_1) 
:=
&   U P \Big( \frac{U^2}{2} + C_V \Theta \Big) 
         + U P R \Theta 
\\
&
+ \frac{d-2}{r_1} \int_0^{r_1}
   \Big( U P \Big( \frac{U^2}{2} + C_V \Theta \Big ) 
       + U P R \Theta 
   \Big)
   dr_2
\\
&
 + \Big( \frac{\kappa}{C_V} - (2\mu + \lambda)\Big)
    \Big( \frac{(U^2)'}{2} + \frac{d-2}{2r_1} U^2 \Big)
  - \lambda (d-1) 
    \Big( \frac{U^2}{r_1}
      + \frac{d-2}{r_1} \int_0^{r_1} \frac{U^2}{r_2}  dr_2 \Big),
\end{split}
\end{equation}
and 
$$
V(r) 
:= \int_0^r \frac{U' + \frac{d-1}{\tilde r} U}
                {\frac{1}{2} \tilde r  - U} d\tilde r ,
\quad 
 W(r) := \frac{1}{2 (2\mu + \lambda)} \int_0^r \tilde r P  (\tilde r)  d\tilde r, 
\quad 
Z (r) := \frac{C_V}{2\kappa} \int_0^r \tilde r  P(\tilde r) d\tilde r.
$$

\noindent 
{\bf Proof of the formulas \eqref{0728-1}, \eqref{0728-2} and 
\eqref{0728-3}. } 
It is easy to check that the first equation in \eqref{eq:ODEs} 
is equivalent to 
$$
\frac{P '}{P} 
= \frac{U' + \frac{d-1}{r} U}{\frac{1}{2}r - U},
$$
which proves \eqref{0728-1} by integrating.

The second equation in \eqref{eq:ODEs} is rewritten as 
$$
-\frac{1}{2} (r P U)' 
+ (P U^2)' + \frac{d-1}{r} P U^2 
- (2\mu + \lambda) 
  \Big( U' + \frac{d-1}{r} U \Big) ' + (P R \Theta)' = 0.
$$
Integrating implies that 
$$
-\frac{1}{2} r P U + P U^2 
 + \int_0^r \frac{d-1}{r_2} P U^2 dr_2
 - (2\mu + \lambda) \Big( U' + \frac{d-1}{r}U \Big)
 + P R \Theta - P_0 R \Theta_0 = 0,
$$
and multiplying by $r^{d-1} e^{W(r)}$ yields that 
$$
- (2\mu + \lambda) 
\Big( r^{d-1} e^{W(r)} U \Big) ' 
+ r^{d-1} e^{W(r)}
 \Big( P U^2 + \int_0^r \frac{d-1}{r_2} P U^2 dr_2
      + P R \Theta - P_0 R \Theta_0 
 \Big)
      = 0.
$$
Hence, we obtain \eqref{0728-2} by integrating the above equation 
and multiplying by $r^{-d+1}e^{-W(r)}$. 

Finally, we consider the third equation in \eqref{eq:ODEs}. 
By multiplying by $r$ and similarly to the argument for 
the second equation, 
we get that 
\begin{equation}\notag 
\begin{split}
& 
-\frac{1}{2} 
 \Big( r^2 P \Big( \frac{U^2}{2} + C_V \Theta \Big) \Big) '
\\
&+ { r^{2-d} }
 \Big( r^{d-1} 
      \Big( U P \Big( \frac{U^2}{2} + C_V \Theta 
                  \Big) 
          + U P R \Theta 
      \Big) 
 \Big) '
 - \kappa r^{-d+2} (r^{d-1} \Theta')'
\\
= 
& r^{-d+2} 
 \Big(  (2\mu + \lambda) (r^{d-1} UU')'
       + \lambda (d-1) ({r^{d-2}} U^2)' 
 \Big) .
\end{split}
\end{equation}
Integrating the above and performing integrations by parts give
\begin{equation}\notag 
\begin{split}
& -\frac{1}{2} r^2 P \Big( \frac{U^2}{2} + C_V \Theta \Big)
\\
& 
+r \Big( U P \Big( \frac{U^2}{2} + C_V \Theta \Big) 
         + U P R \Theta 
   \Big)
+ (d-2) \int_0^r 
   \Big( U P \Big( \frac{U^2}{2} + C_V \Theta \Big ) 
       + U P R \Theta 
   \Big)
   dr_2
\\
&
 - \kappa \Big( r \Theta ' + (d-2) \Big(\Theta  (r) - \Theta_0 
 \Big) \Big)
\\
= 
& (2\mu + \lambda) 
  \Big( r UU' + (d-2) \int_0^r UU' d r_2 \Big)
  + \lambda (d-1) 
   \Big( U^2 + (d-2) \int_0^r \frac{U^2}{r_2}  dr_2 \Big) .
\end{split}
\end{equation}
Dividing by $r$, and regarding this formula as an equation on the energy $U^2/2 + C_V \Theta$, we get
\begin{equation}\notag 
\begin{split}
& -\frac{1}{2} r P \Big( \frac{U^2}{2} + C_V \Theta \Big)
\\
& 
+ U P \Big( \frac{U^2}{2} + C_V \Theta \Big) 
         + U P R \Theta 
+ \frac{d-2}{r} \int_0^r 
   \Big( U P \Big( \frac{U^2}{2} + C_V \Theta \Big ) 
       + U P R \Theta 
   \Big)
   dr_2
\\
&
 - \frac{\kappa}{C_V} 
   \Big( \Big( \frac{U^2}{2} + C_V\Theta \Big) ' 
       + \frac{d-2}{r} \Big( \frac{U^2}{2} + C_V \Theta \Big) 
   \Big)
\\
& + \frac{\kappa}{C_V}
    \Big( \frac{(U^2)'}{2} + \frac{d-2}{2r} U^2 \Big)
  + \kappa \frac{d-2}{r} \Theta_0
\\
= 
& (2\mu + \lambda) 
  \Big( UU' + \frac{d-2}{2r} U^2 \Big)
  + \lambda (d-1) 
   \Big( \frac{U^2}{r}
      + \frac{d-2}{r} \int_0^r \frac{U^2}{r_2}  dr_2 \Big) .
\end{split}
\end{equation}
Multiplying by $r^{d-2} e^{Z(r)}$ gives
\begin{equation}\label{0728-4}
\begin{split}
& \frac{\kappa}{C_V} 
  \Big( r^{d-2} e^{Z(r)} \Big( \frac{U^2}{2} + C_V \Theta \Big)
  \Big) '
\\
= 
& 
  r^{d-2} e^{Z(r)} 
 \Big\{ U P \Big( \frac{U^2}{2} + C_V \Theta \Big) 
         + U P R \Theta 
+ \frac{d-2}{r} \int_0^r 
   \Big( U P \Big( \frac{U^2}{2} + C_V \Theta \Big ) 
       + U P R \Theta 
   \Big)
   dr_2
\\
& + \Big( \frac{\kappa}{C_V} - (2\mu + \lambda)\Big)
    \Big( \frac{(U^2)'}{2} + \frac{d-2}{2r} U^2 \Big)
  - \lambda (d-1) 
    \Big( \frac{U^2}{r}
      + \frac{d-2}{r} \int_0^r \frac{U^2}{r_2}  dr_2 \Big)
\\
&  + \kappa \frac{d-2}{r} \Theta_0 \Big\} 
\\
=: 
&
r^{d-2} e^{Z(r)}  
   \Big( F_\Theta (r)  + \kappa \frac{d-2}{r} \Theta_0 \Big) ,
 \end{split}
\end{equation}
where $F_\Theta$ is also given in the formula \eqref{0728-3}.  Integrating the above leads to
\begin{equation}\notag 
\begin{split}
& \frac{\kappa}{C_V} 
  r^{d-2} e^{Z(r)} \Big( \frac{U^2}{2} + C_V \Theta \Big)
= 
\int_0^r 
  r_1^{d-2} e^{Z(r_1)} 
 \Big( 
 F_\Theta (r_1)  
  + \kappa \frac{d-2}{r_1} \Theta_0 \Big)  \, dr_1 , 
\end{split}
\end{equation}
which proves \eqref{0728-3}.
\hfill $\Box$

\subsection{Integro-differential equations for cavitating solutions}
Let us consider the vacuum case at $x = 0$.  
Supposing the condition at $r = 0$ that for 
$0 < \alpha <  1/2$
\[
P(0) = 0, \qquad 
U(0) = 0, \qquad U'(0) = \alpha , \qquad 
\Theta (0) = \Theta_0 > 0 , \qquad  \Theta '(0) = 0
\]
and $P, U, \Theta$ are $C^1$, we get the following integral equations
\begin{align}\label{0728-1vacuum}
P (r) = & 
e^{V(r) - V(\delta)} P (\delta) 
\text{ for given } \delta > 0,
\\ \label{0728-2vacuum}
U(r) 
=
& \frac{r^{-d+1}}{2\mu + \lambda} 
\int_0^r 
  r_1 ^{d-1} e^{-W(r) + W(r_1)} \widetilde F_U(r_1) \, dr_1,
\\ \notag 
\Theta (r) 
= 
&
\displaystyle 
(d-2)r^{-d+2}
\int_0^r r_1 ^{d-3} e^{-Z(r) + Z (r_1)} \, dr_1 \Theta_0 - \frac{U^2}{2C_V}
\\ \label{0728-3vacuum} 
& 
+ \frac{r^{-d+2}}{\kappa} \int_0^r 
  r_1^{d-2} e^{-Z(r) + Z(r_1)} F_\Theta (r_1)
 \, dr_1 ,
\end{align}
where $F_\Theta, V(r), W(r), Z(r)$ are same as in \eqref{0728-1}, \eqref{0728-2}, \eqref{0728-3} and 
\begin{equation} \notag 
\begin{split}
\widetilde F_U(r_1) 
:=
&      P U^2 
       + \int_0^{r_1} \frac{d-1}{r_2} P U^2 \, dr_2 
       + P R \Theta + d (2\mu + \lambda) \alpha .  
\end{split}
\end{equation}
On the density $P(r)$, 
the exponent $V(r)$ itself diverges 
because of $V'(r) \simeq \frac{d\alpha}{\frac{1}{2} -\alpha}\cdot \frac{1}{r}$, 
but we always regard 
$V(r) - V(\delta)$ as an integral from $\delta$ to $r$.

\section{Smooth self-similar solutions}

\subsection{Main result}

\begin{notation}
In this section, we consider $R, \mu, \lambda, C_V, P_0$ as  fixed positive constants - the meaning of $P_0$ will soon be explained. We denote $C$ a constant which depends on $(R,\mu,\lambda,C_V,P_0) \in (0,\infty)^5$; the implicit constant in the notations $\lesssim$ and $O(\cdot)$ has the same properties.
\end{notation}

\begin{thm} \label{theo:2}
For fixed $R, \mu, \lambda, C_V, P_0>0$, if $\Theta_0$ is sufficiently small, there exists a unique continuously differentiable function $(P,U,\Theta) \in C^1 ([0,\infty))^3  $ solving~\eqref{eq:ODEs} such that, for $r$ small,
\begin{align*}
& P(0) = P_0, \\
& U(r) = O(r^2), \qquad U'(r) = O(r), \\
& \Theta(r) = \Theta_0, \qquad \Theta'(r) = O(r).
\end{align*}
It satisfies the bounds
\begin{align*}
|U(r)| \lesssim \frac{\Theta_0 r^2}{(1+r)^3}, \qquad |U'(r)| \lesssim \frac{\Theta_0 r}{(1+r)^3}, \\
|\Theta(r)| \lesssim \frac{\Theta_0}{(1+r)^2}, \qquad |\Theta'(r)| \lesssim \frac{\Theta_0 r}{(1+r)^2}.
\end{align*}
Furthermore, there exists $P_\infty>0$, $U_\infty < 0$, $\Theta_\infty>0$ such that
\begin{align*}
& P(r) = P_\infty + O \left( \frac{1}{r^2} \right), \\
& U(r) = \frac{U_\infty}{r} + O \left( \frac{1}{r^3} \right) ,\\
& \Theta(r) = \frac{\Theta_\infty}{r^2} + O\left( \frac{1}{r^4} \right).
\end{align*}
Finally,
\begin{align*}
& U_{\infty} = - 2R \Theta_0 + O (\Theta_0^2), \\
& \Theta_\infty = \frac{2(d-2) \kappa}{C_V P_0} \Theta_0 + 
O (\Theta_0^2 ).
\end{align*}
\end{thm}

\subsection{Main steps of the proof}

Let
\begin{equation}\notag 
\begin{split}
\| (U,\Theta) \|^{\delta}
:= 
&
\sup _{0<r<\delta} \left[ r^{-2} |U(r)| 
+ r^{-1} |U'(r)| + |\Theta(r)| + r^{-1} |\Theta' (r)|\right].
\end{split}
\end{equation}

To prove the existence of a local solution, let $\Psi$ be the map which to $(U,\Theta)$ associates the right-hand side of~\eqref{0728-2} and~\eqref{0728-3}:
$$
\Psi \; : \; (U,\Theta) \mapsto (\operatorname{RHS}\eqref{0728-2},\operatorname{RHS}\eqref{0728-3}) . 
$$
Define
$$
E^\delta = \{ (U,\Theta) \in \mathcal{C}^1(0,\delta) \; \; \mbox{such that} \;\; \Theta(0) = \Theta_0 \;\; \mbox{and} \;\; \| (U,\Theta) \|^\delta < \infty \}.
$$
Equipped with $\| \cdot \|^\delta$, it is a Banach space (an affine Banach space to be precise).

\begin{lem}
\label{lempsi}
There exists $C_0 >0$ such that: setting $\epsilon = C_0 \Theta_0$, if $\delta$ and $\epsilon$ are sufficiently small, then $\Psi$ is a contraction on the ball $B_{E^\delta} ((0,\Theta_0),\epsilon)$.
\end{lem}

By the Banach fixed point theorem, this lemma gives the local existence (close to zero) of solutions. They can then be prolonged for $r>0$ as long as $U,U',\frac{1}{U - \frac{1}{2}r},\Theta,\Theta'$ are bounded. 
Let $[0,R)$ be the largest interval on which $P,U,\Theta$ are well-defined. In other words, either $R = \infty$, or
\begin{equation}
\label{defR1}
\lim_{r \to R} |U(r)| + |U'(r)| + \frac{1}{|U - \frac{1}{2}r|}  + |\Theta(r)| + |\Theta'(r)| = \infty.
\end{equation}

Define then, for constants $M_1$ and $M_2$ which are much smaller than $1$,
\begin{align*}
Z(s) = \sup_{0<r<s} & \frac{1}{M_1} r^{-2} (1+r)^3 |U(r)| 
+ \frac{1}{M_1}r^{-1}(1+r)^3 \left|U'(r) + \frac{d-1}{r} U(r) \right| \\
& \qquad + \frac{1}{M_2} (1+r)^2 |\Theta(r)| + \frac{1}{M_2}r^{-1} (1+r)^2 |\Theta' (r)|.
\end{align*}

\begin{lem}
\label{lemZ}
If $Z(r) \leq 1$, then
$$
Z(r) \lesssim M_1 + \frac{M_2}{M_1} + \frac{\Theta_0}{M_1}+ \frac{\Theta_0}{M_2} + \frac{M_1^2}{M_2}.
$$
\end{lem}

We now choose the constants $M_1$ and $M_2$ such that
\begin{equation}
\label{defTMM}
M_1^2 + M_1 M_2 \ll \Theta_0  \ll M_2 \ll M_1 \ll 1.
\end{equation}
(where the notation $A \ll B$ means $A \leq cB$, for a constant $c$ depending on the parameters of the problem ($R$, $\mu$, $\lambda$, $C_V$, $P_0$), which is chosen sufficiently small 
 so 
that all arguments in the following apply). For instance, to achieve the above, we could choose
$$
M_2 = A \Theta_0, \qquad M_1 = A^2 \Theta_0.
$$
Choosing $A$ sufficiently large ensures that $\Theta_0 \ll M_2 \ll M_1$; choosing then $\Theta_0 \ll \frac{1}{A^4}$ ensures that $M_1 \ll 1$ and $M_1^2 + M_1 M_2 \ll \Theta_0$.
Let 
$$
\widetilde{R} = \sup \{ r \; \mbox{such that} \; Z(r) \leq 1\}.
$$
By Lemma~\ref{lempsi}, $\widetilde R \geq \delta$. Argue by contradiction and assume that $\widetilde{R}$ is finite. Then, by~\eqref{defR1}, $\widetilde{R} < R$. Furthermore, by Lemma~\ref{lemZ} and the choice~\eqref{defTMM}, $Z(\widetilde{R}) < 1$. Then $(U,\Theta)$ can be prolonged over a short time interval where $Z < 1$. This contradicts the definition of $\widetilde{R}$, and gives the desired result: $R = \widetilde{R} = \infty$.

There remains to prove the asymptotic behavior of $P,U,\Theta$. This is achieved in the following lemma.
\begin{lem} \label{asymptoticbehav}
There exists $P_\infty>0$, $U_\infty \in \mathbb{R}$,  $\Theta_\infty>0$ such that
\begin{align*}
P(r) = P_\infty + O \left( \frac{1}{r^2} \right), \; U(r) = \frac{U_\infty}{r} + O \left( \frac{1}{r^3} \right), \;\; \Theta(r) = \frac{\Theta_\infty}{r^2} + O\left( \frac{1}{r^4} \right).
\end{align*}
Furthermore, $\Theta(r)>0$ for any $r$ and
$$
U_\infty = -2R \Theta_0 + O(\Theta_0^2) \qquad \Theta_\infty = \frac{2(d-2) \kappa}{C_V P_0} \Theta_0 + O(\Theta_0^2).
$$
\end{lem}

\subsection{Local existence: proof of Lemma~\ref{lempsi}}
We assume here that $(U_1,\Theta_1), (U_2,\Theta_2) \in B_{E^\delta}(0,\epsilon)$, and let
$$
D = \| (U_1,\Theta_1) - (U_2,\Theta_2) \|^\delta.
$$
We aim at proving that
\begin{equation}
\label{eqD}
\| \Psi (U_1,\Theta_1) - \Psi (U_2,\Theta_2) \|^\delta \lesssim (\epsilon + \delta) D.
\end{equation}
This proves that $\Psi$ acts as a contraction for $\epsilon$, $\delta$ sufficiently small. Furthermore, this proves that $\Psi$ stabilizes $B_{E^\delta}((0,\Theta_0),\epsilon)$. Indeed, choosing $C_0$ sufficiently large, this follows from~\eqref{eqD} together with the observation that $(0,\widehat \Theta) := \Psi(0,\Theta_0)$ satisfies
$$
\widehat{\Theta}(r) = \Theta_0(1 + O(r^2)) \qquad \mbox{and} \qquad |\widehat{\Theta}'(r)| \lesssim \Theta_0 r.
$$ 
There remains to prove~\eqref{eqD}!

\bigskip

With the notation used in~\eqref{0728-2} and ~\eqref{0728-3}, it appears first that, as soon as $\|(U,\Theta)\|^\delta < \epsilon$,
\begin{align*}
& |V| \lesssim 1, \quad |W| + |Z| \lesssim r^2\\
& |V'| \lesssim \epsilon, \quad |W'| \lesssim r, \quad |Z'| \lesssim r \\
& |P| \lesssim 1, \quad |P'| \lesssim \epsilon \\
&|F_U(r)| \lesssim \epsilon r\\
&|F_\Theta(r)| \lesssim \epsilon^2 r^2.
\end{align*}

We now turn to the difference between two solutions. We denote in the following $(\widetilde{U}_1,\widetilde{\Theta}_1) = \Psi (U_1,\Theta_1)$, etc... By direct inspection,
\begin{align*}
|V_1 - V_2| + |P_1 - P_2| \lesssim Dr, \qquad |W_1 - W_2| + |Z_1 - Z_2| \lesssim D r^3.
\end{align*}
It follows that
\begin{align*}
& |\Theta_1(r) - \Theta_2(r)| \leq \int_0^r |\Theta_1' - \Theta_2'|\,ds \lesssim \int_0^r Ds\,ds \lesssim Dr^2 \\
& |F_{U_1}(r) - F_{U_2}(r)| \lesssim D (\epsilon + \delta)  r \\
& |F_{\Theta_1}(r) - F_{\Theta_2}(r)| \lesssim D \epsilon r^2.
\end{align*}
Therefore,
\begin{align*}
|\widetilde U_1(r) - \widetilde U_2(r)| \lesssim & r^{1-d} \int_0^r r_1^{d-1} |e^{W_1(r_1)  - W_1(r)} - e^{W_2(r_1)  - W_2(r)}| F_{U_1}(r_1) \, dr_1 \\
& \qquad + r^{1-d} \int_0^r r_1^{d-1} e^{W_2(r_1)  - W_2(r_1)} |F_{U_1}(r) - F_{U_2}(r_1)| \,dr_1 \\
& \lesssim r^{1-d} \int_0^r r_1^{d-1} \left[ D (r^3 + r_1^3) \epsilon r  + D(\epsilon + \delta) r_1 \right]\,dr_1 \lesssim D(\epsilon + \delta) r^2.
\end{align*}
Arguing similarly,
$$
|\widetilde U_1'(r) - \widetilde U_2'(r)| \lesssim D(\epsilon + \delta) r
$$
and finally
$$
|\widetilde \Theta_1(r) - \widetilde \Theta_2(r)| \lesssim D\epsilon r^3 \qquad \mbox{and} \qquad |\widetilde \Theta_1'(r) - \widetilde \Theta_2'(r)| \lesssim D\epsilon r^2,
$$
which concludes the proof.

\subsection{Global existence: proof of Lemma~\ref{lemZ}}
We assume here that we have a solution defined on $[0,R_0]$  for $R_0 > 0$, such that $Z(r) \leq 1 $ for all $ r \in [0,R_0]$ which implies
\begin{align}\label{0914-1}
& |U(r)| \leq  M_1 r^2 (1+r)^{-3}, \quad |U'(r)| \lesssim M_1 r (1+r)^{-3},\\
\label{0914-2}
&  |\Theta(r)| \lesssim M_2 (1+r)^{-2}, \quad |\Theta'(r)| \lesssim M_2 r (1+r)^{-2}.
\end{align}

\bigskip
\noindent 
\underline{Estimate of $P$. } 
It follows from the definition of $V$ that 
\begin{equation}\label{0728-7}
|V(r)|  
\leq \int_0^r 
\frac{\frac{s}{(1+s)^3}M_1 + \frac{d-1}{s}\frac{s^2}{(1+s)^3}M_1}
     {\frac{1}{2}s - \frac{s^2}{(1+s)^3}M_1} 
     ds
\leq \int_0^r\frac{d\frac{s}{(1+s)^3}M_1}{s (\frac{1}{2}-M_1)} 
     ds
\leq \frac{dM_1}{1-2M_1} =: M_0 .
\end{equation}
Hence we obtain that 
\begin{gather}\label{0728-5}
e^{-M_0} P_0 \leq P (r) \leq e^{M_0} P_0, 
\quad r > 0, 
\\ \label{0728-5(2)}
|P ' (r)| 
\leq  |V'(r)| P (r)
\leq \frac{d\frac{r}{(1+r)^3}M_1}{r (\frac{1}{2}-M_1)}  e^{M_0}P_0
\leq \frac{2M_0}{(1+r)^3} e^{M_0} P_0, 
\quad r > 0 .
\end{gather}

\bigskip

\noindent 
\underline{Estimate of $U$. }
We have from \eqref{0728-5} that 
\begin{equation}\label{0913-1}
- e^{M_0} P_0 (r^2 -r_1^2)
\lesssim -W(r) + W(r_1) 
\lesssim - e^{-M_0} P_0(r^2 -r_1^2), 
\quad r \geq r_1 .
\end{equation}
Furthermore, by the above bound on $P$ and $P'$,
$$
|R P(r) \Theta(r) - R P_0 \Theta_0| \lesssim M_2 \frac{r}{1+r}.
$$
As a consequence,
\begin{equation}
\label{boundfu}
F_U(r) \lesssim (M_1^2 + M_2) \frac{r}{1+r} \qquad \mbox{and} \qquad |F_U'(r)| \lesssim \frac{M_1^2 + M_2}{1+r}.
\end{equation}
Therefore,
$$
|U(r)| \lesssim  (M_1^2 + M_2) r^{1-d} \int_0^r r_1^{d-1} e^{-C (r^2-r_1^2)}  \frac{r_1}{1+r_1}\,dr_1.
$$
Using the inequality (for $C>0$)
\begin{equation}
\label{basicinequality}
\int_1^r s^{\alpha} e^{-C(r^2-s^2)} \,ds \lesssim r^{\alpha-1},
\end{equation}
we get that
$$
|U(r)| \lesssim (M_1^2 + M_2) \frac{r^2}{(1+ r)^3}.
$$
As for the derivative of $U$, 
we can get the required estimates for small $r $ easily. 
In fact, we directly differentiate to estimate 
$$
|U'(r)| \lesssim ( M_1^2 + M_2 )  r , 
 \quad 0 < r \leq 1 .
$$
In order to deal with $r>1$, we write 
\begin{equation}\label{0728-10}
\begin{split}
U'(r) + \frac{d-1}{r}  U(r)
= \frac{F_ U (r)}{2\mu + \lambda}
 - W'(r) U(r) , 
\end{split}
\end{equation}
In order to see that the right-hand side is decaying, we must take advantage of a cancellation between the two terms. It becomes apparent after integrating by parts:
\begin{equation}\notag 
\begin{split}
-W'(r) U(r) 
= 
& 
 -W'(r) \frac{r^{-d+1}}{2\mu + \lambda} 
 \int_0^r r_1 ^{d-1} 
   \frac{\partial_{r_1} e^{-W(r)+W(r_1)} }{W'(r_1)} 
   F_ U(r_1) \, dr_1
\\
=
& -\frac{ F_ U(r)}{2\mu + \lambda} 
  +W'(r) \frac{r^{-d+1}}{2\mu+\lambda} 
   \int_0^r e^{-W(r)+W(r_1)} 
     \partial_{r_1} 
     \Big( \frac{r_1^{d-1} F_ U(r_1)}{W'(r_1)} \Big) dr_1 ,
\end{split}
\end{equation}
so that
\begin{equation}
\label{doublestar}
U'(r) + \frac{d-1}{r} U(r) = W'(r) \frac{r^{-d+1}}{2\mu + \lambda}  \int_0^r e^{-W(r)+W(r_1)} 
     \partial_{r_1} 
     \Big( \frac{r_1^{d-1} F_ U(r_1)}{W'(r_1)} \Big) dr_1. 
\end{equation}
By \eqref{0728-5}, \eqref{0728-5(2)} and~\eqref{boundfu},
\begin{equation}\label{0913-2}
\begin{split} 
\Big|
     \partial_{r_1} 
     \Big( \frac{r_1^{d-1} F_U (r_1)}{W'(r_1)} \Big)
\Big| 
= 
& 
{2(2\mu + \lambda)} \Big|  
  \frac{ (d-1)r_1^{d-2} F_ U(r_1) +r_1 ^{d-1} F_  U'(r_1)}
  {r_1 P (r_1)}
\\
& \quad + r_1 ^{d-1} F_ U(r_1)
  \Big( -r_1^{-2} P (r_1)^{-1} 
       - r_1 ^{-1} P ' (r_1) P (r_1)^{-2}
  \Big)
\Big|
\\
\lesssim
& 
r_1^{d-3}(M_1^2 + M_2)
\end{split}
\end{equation}
Therefore, for $r>1$,
\begin{equation}\notag 
\begin{split}
\Big| U'(r) + \frac{d-1}{r} U(r) \Big| 
\lesssim
& (M_1^2 + M_2) r r^{-d+1} 
\int_0^r e^{-C (r^2 - r_1 ^2)} r_1 ^{d-3}
\, dr_1 \lesssim (M_1^2 + M_2) r^{-2},
\end{split}
\end{equation}
where the last inequality follows from~\eqref{basicinequality}.
Summarizing, we obtain the estimate
\begin{equation}\label{0728-11}
\begin{split}
\left|U'(r) + \frac{d-1}{r} U(r) \right|
\lesssim (M_1^2 + M_2) \frac{r}{(1+r)^3}.
\end{split}
\end{equation}

\bigskip

\noindent
\underline{Estimate of $\Theta$}. 
Observe that
$$
|F_\Theta(s)| \lesssim (M_1^2 + M_1 M_2) \frac{s}{(1+s)^2}
$$
while if $0 \leq r_1 \leq r$,
$$
 -Z(r) + Z(r_1) \leq C (r_1^2 - r^2).
$$
It follows then from~\eqref{basicinequality} that
\begin{equation*}
|\Theta(r)| \lesssim  \frac{\Theta_0 + M_1^2 + M_1 M_2}{(1+r)^2} . 
\end{equation*}
Finally, noticing that 
\begin{equation}\label{0914-5}
(d-2)r^{-d+2} \int_0^r r_1 ^{d-3} e^{-Z(r) + Z (r_1)} \, dr_1 = 1 + (d-2)r^{-d+2}\int_0^r r_1 ^{d-3} \left[ e^{-Z(r) + Z (r_1)} - 1 \right]\, dr_1, 
\end{equation}
one can prove that
$$
|\Theta'(r)| \lesssim ( \Theta_0 + M_1^2 + M_1 M_2) \frac{r}{(1+r)^2}.
$$

\subsection{Asymptotic behavior: proof of Lemma~\ref{asymptoticbehav}}
\noindent \underline{Asymptotic behavior of $P$.} 
It follows from the estimate \eqref{0728-5(2)} that 
$$
\int_0^\infty |P'(r)| dr < \infty \qquad \mbox{and} \qquad P(r) = P_\infty + O \left( \frac{1}{r^2} \right) 
\text{ with } P_\infty = \int_0^\infty P'(r) dr .
$$

\medskip

\noindent \underline{Asymptotic behavior of $U$.} 
The terms $PU^2 + PR\Theta $, which are parts of $F_U$, can be estimated by
\begin{equation}\notag 
\begin{split}
\frac{r^{-d+1}}{2\mu + \lambda} 
\int_0^r r_1 ^{d-1} e^{-W(r) + W(r_1)} 
   \big| PU^2 + PR\Theta \big| d\tau 
\lesssim 
&  r^{-d+1} 
   \int_0^r r_1 ^{d-1} e^{-c (r^2 - r_1^2)} 
     \frac{1}{(1+r_1)^2} \,dr_1
\lesssim 
  r^{-3}.
\end{split}
\end{equation}
The other terms can be written after integration by parts
\begin{equation}\label{0809-2} 
\begin{split}
& 
\frac{r^{-d+1}}{2\mu + \lambda} 
\int_0^r r_1 ^{d-1} 
   \frac{2(2\mu + \lambda)}{r_1 P(r_1)}
   \Big( \partial_{r_1}e^{-W(r) + W(r_1)} \Big)
   \Big( \int_0^{r_1} \frac{d-1}{r_2} PU^2 dr_2 - P_0 R \Theta_0 
   \Big) d\tau 
\\
& \qquad =
 \frac{2 r^{-1}}{P(r)} 
   \Big( \int_0^{r} \frac{d-1}{r_2} PU^2 dr_2 - P_0 R \Theta_0
   \Big) 
   + A(r),
\end{split}
\end{equation}
where 
\begin{equation}\notag 
\begin{split}
A(r)=
2 r^{-d+1}
\int_0^r 
   e^{-W(r) + W(r_1)} 
   \partial_{r_1} 
   \Big( \frac{r_1 ^{d-2} }{ P(r_1)}
   \Big( \int_0^{r_1} \frac{d-1}{r_2} PU^2 dr_2 - P_0 R \Theta_0 
   \Big) 
   \Big) 
   d\tau.
\end{split}
\end{equation}
Noting that
$$
\int_0^{\infty} \Big|\frac{d-1}{r_2} PU^2 \Big| dr_2 < \infty\quad \mbox{and} \quad \int_0^r \frac{d-1}{r_2} PU^2 \, dr_2 = \int_0^\infty \frac{d-1}{r_2} PU^2 \, dr_2 + O \left( \frac{1}{r^2} \right),
$$
we see that, for a constant $U_\infty$,
$$
\frac{2 r^{-1}}{P(r)} 
   \Big( \int_0^{r_1} \frac{d-1}{r_2} PU^2 dr_2 - P_0 R \Theta_0
   \Big) 
 = \frac{U_\infty}{r} + O \left( \frac{1}{r^3} \right).
$$
Finally, it follows from~\eqref{basicinequality} that $A(r) \lesssim r^{-3}$. 

\medskip

\noindent \underline{Asymptotic behavior of $\Theta$.} 
Similarly to the above argument for $U$, 
we deduce that $\Theta$ can be written
\begin{align*}
 \Theta(r) & =
(d-2)r^{-d+2}
\int_0^r r_1 ^{d-3} e^{-Z(r) + Z (r_1)} \, dr_1 \Theta_0
- \frac{U^2}{2C_V} \\
& \qquad \qquad \qquad \qquad \qquad \qquad
 + \frac{r^{-d+2}}{\kappa} \int_0^r 
  r_1^{d-2} e^{-Z(r) + Z(r_1)} \widetilde F_\Theta (r_1)
 \, dr_1 + O\left( \frac{1}{r^4} \right) , 
\end{align*}
where 
$$
\widetilde F_\Theta (r_1) 
:=
\frac{d-2}{r_1} \int_0^{r_1}
   \Big( U P \Big( \frac{U^2}{2} + C_V \Theta \Big ) 
       + U P R \Theta 
   \Big)
   dr_2
  - \frac{ \lambda (d-1)(d-2)}{r_1} 
   \int_0^{r_1} \frac{U^2}{r_2}  dr_2 .
$$
The integration by parts allows us to handle these terms 
in a similar way to that for $U $ and we obtain 
that there exists $C_2$ such that 
$$
\widetilde \Theta (r) = C_2 r^{-2} + O (r^{-4}) 
\quad \text{as } r \to \infty .
$$

\bigskip 

\noindent 
\underline{Positivity of temperature for all $r > 0$}. 
It is easy to check that there exists $c > 0$ such that 
\begin{equation}\label{0906-1}
(d-2) r^{-d+2} \int_0^r 
 r_1 ^{d-3} e^{-Z(r) + Z(r_1)} dr_1 \Theta _0 
 \geq c \Theta _0 (1+r)^{-2}  
 \quad \text{for any } r > 0. 
\end{equation}
On the other hand, we can estimate the others as 
\begin{equation}\label{0906-2}
\Big| 
\Theta (r) - 
(d-2) r^{-d+2} \int_0^r 
 r_1 ^{d-3} e^{-Z(r) + Z(r_1)} dr_1 \Theta _0 
\Big| 
\lesssim (M_1^2 + M_1 M_2) (1+r) ^{-2} \ll  \Theta_0 (1+r) ^{-2}, 
\end{equation}
which proves the positivity of $\Theta$ provided that 
$\Theta _0 \ll 1$.

\bigskip

\noindent \underline{Expansion of $\Theta_\infty$ and $U_\infty$.} Observe first that
$$
\Theta(r) = (d-2) r^{-d+2} \int_0^r r_1^{d-3} e^{-Z(r) + Z(r_1)} \,dr_1 \Theta_0 + O \left( \frac{\Theta_0^2}{ r^2} \right) + O \left( \frac{1}{r^3} \right).
$$
Since $P_\infty = P_0+ O(\Theta_0)$, we get $Z(r) = \frac{C_V P_0}{4\kappa} r^2 + O(\Theta_0 r^2)$, and therefore, as $r \to \infty$,
$$
(d-2) r^{-d+2} \int_0^r r_1^{d-3} e^{-Z(r) + Z(r_1)} \,dr_1 \Theta_0 = \left[ \frac{2(d-2)\kappa}{C_V P_0} \Theta_0 + O (\Theta_0 ^2)
 \right] \frac{1}{r^2}.
$$
This means that
$$
\Theta_\infty = \frac{2(d-2)\kappa}{C_V P_0} \Theta_0 + O (\Theta_0^2).
$$
Arguing similarly, for $U$, one finds $U_\infty = - 2R \Theta_0 +O(\Theta_0^2)$.

\section{Cavitating self-similar solutions}

\subsection{Main result}

\begin{notation}
In this section, we consider $R, \mu, \lambda, C_V$ as  fixed positive constants. We denote $C$ a constant which depends on $(R,\mu,\lambda,C_V) \in (0,\infty)^4$; the implicit constant in the notations $\lesssim$ and $O(\cdot)$ has the same properties.

Furthermore, for given quantities $A$ and $B$, we denote $A \ll B$ to mean that $A \leq c B$ for a constant $c = c(R,\mu,\lambda,C_V)$ which is sufficiently small so that all the needed arguments apply.
\end{notation}

\begin{thm} \label{theo:3}
There exists a constant $\epsilon > 0$ such that: if
\begin{equation}
\label{smallnesscondition}
P_\delta + \delta + \Theta_0 + \alpha + \frac{P_\delta \Theta_0}{\alpha} + \frac{\alpha^2}{P_\delta  \Theta_0} + \alpha \log \frac{1}{P_\delta \delta^2} < \epsilon,
\end{equation}
then there exists a solution $(P,U,\Theta) \in C^1 ([0,\infty)) \times C^1([0,\infty)) \times C^{\frac{2d\alpha}{1-2\alpha}}$ solving~\eqref{eq:ODEs} such that $P(\delta) = P_\delta$, $U(0) = 0$, $U'(0) = \alpha$, $\Theta(0) = \Theta_0$, $\Theta'(0) = 0$. For $r$ small, it satisfies
\begin{align*}
& P(r) = P_\delta \left( \frac{r}{\delta} \right)^{\frac{2d\alpha}{1-2\alpha}} + O \left( \frac{r}{\delta} \right)^{\frac{2d\alpha}{1-2\alpha}+1+d\alpha},\\
& U(r) = \alpha r + O(r^{1 + \frac{2d\alpha}{1-2\alpha}}), \\
& \Theta(r) = \Theta_0 + O(r^2).
\end{align*}
It also satisfies the global bounds
\begin{align*}
& |P(r)| \lesssim P_\delta \min\left[1 , \left( \frac{r}{\delta} \right)^{\frac{2d\alpha}{1-2\alpha}} \right], \\
& |U(r)| \lesssim \frac{\alpha r}{(1+ \sqrt{P_\delta} r)^2}, \qquad |U'(r)| \lesssim \frac{\alpha}{(1+\sqrt{P_\delta} r)^2}, \\
& |\Theta(r)| \lesssim \frac{1}{(1+\sqrt{P_\delta}r)^2}, \qquad |\Theta'(r)| \lesssim \frac{\sqrt{P_\delta} r}{(1+\sqrt{P_\delta}r)^2}.
\end{align*}
Furthermore, there exists $P_\infty>0$, $U_\infty < 0$, $\Theta_\infty>0$ such that
\begin{align*}
& P(r) = P_\infty + O \left( \frac{1}{r^2} \right), \\
& U(r) = \frac{U_\infty}{r} + O \left( \frac{1}{r^3} \right), \\
& \Theta(r) = \frac{\Theta_\infty}{r^2} + O\left( \frac{1}{r^4} \right).
\end{align*}
Finally,
\begin{align*}
& U_{\infty} =  
\frac{2d(2\mu + \lambda)}{P_\delta} \alpha 
    \left[ 1 + O\Big( \alpha \log \frac{1}{P_\delta \delta^2} \Big) \right], \\
& \Theta_\infty =  
\frac{2(d-2)\kappa}{C_V P_\delta} \Theta_0 
\left[ 1+ O\left(\alpha \log \frac{1}{P_\delta \delta ^2}\right) 
 + O\left( \frac{\alpha^2 }{P_\delta \Theta_0} \right) 
 \right]. 
\end{align*}
\end{thm}

\begin{rem} The above result remains true if $\Theta_0$ is not assumed to be small as in~\eqref{smallnesscondition}, but only $O(1)$.
\end{rem}

\begin{rem}
To clarify the meaning of~\eqref{smallnesscondition}, we can scale $P_\delta$ and $\Theta_0$ in terms of $\alpha$
$$
P_\delta = \alpha^{\nu_1}, \qquad \Theta_0 = \alpha^{\nu_2},
$$
while keeping $\delta \ll 1$ fixed. Then~\eqref{smallnesscondition} amounts to requiring that $\alpha \ll 1$, together with
$$
\nu_1 >0, \qquad \nu_2>0, \qquad \mbox{and} \qquad 1< \nu_1 + \nu_2 < 2.
$$
\end{rem}

\subsection{Main steps of the proof}

Let
\begin{equation}\notag 
\begin{split}
\| (U,\Theta) \|^{\delta}
:= 
&
\sup _{0<r<\delta} \left[ r^{-1} |U(r)| 
+ |U'(r)| + |\Theta(r)| + r^{-1} |\Theta' (r)|\right].
\end{split}
\end{equation}

To prove the existence of a local solution, let $\Psi$ be the map which to $(U,\Theta)$ associates the right-hand side of~\eqref{0728-2} and~\eqref{0728-3}:
$$
\Psi \; : \; (U,\Theta) \mapsto (\operatorname{RHS}\eqref{0728-2vacuum},\operatorname{RHS}\eqref{0728-3vacuum}) . 
$$
Define
$$
E^\delta = \{ (U,\Theta) \in \mathcal{C}^1(0,\delta) \; \; \mbox{such that} \;\; \Theta(0) = \Theta_0 \;\; \mbox{and} \;\; \| (U,\Theta) \|^\delta < \infty \}.
$$
This is an affine Banach space.
\begin{lem}
\label{lempsivacuum}
If $\alpha + \delta + P_\delta + \Theta_0 +  \frac{P_\delta \Theta _0}{\alpha} $ is sufficiently small, $\psi$ is a contraction on $B_{E^\delta}((\Theta_0,\alpha r), \alpha/2)$.
\end{lem}

By the Banach fixed point theorem, this lemma gives the local existence (close to zero) of solutions. In order to prolong them, we will argue as in the case of smooth self-similar solutions and define a $Z$ function.

Observing that
\begin{equation}
\begin{split} \label{oriole} 
& r^{-d+1} \int_0^r r_1 ^{d-1} e^{- C P_\delta(r^2 - r_1^2)} dr_1 
\simeq \frac{r}{(1 + \sqrt{P_\delta} r)^2} , 
\\ 
& r^{-d+2} 
\int_0^r r_1 ^{d-3} e^{-CP_\delta (r^2-r_1^2)}dr_1 
\simeq \frac{1}{(1+\sqrt{P_\delta} r)^2}, 
\\
& \left| \partial_r  \left[ r^{-d+2} 
\int_0^r r_1 ^{d-3} e^{-CP_\delta (r^2-r_1^2)}dr_1 \right]
 \right| \lesssim \frac{P_\delta r}{1 + P_\delta r^2} , 
\end{split}
\end{equation}
we are led to defining
\[
\begin{split}
Z(s) 
:= 
& 
\sup _{0 < r < s} 
\frac{1}{M_1}r^{-1}\left(1+ \sqrt{P_\delta}r \right)^2 |U(r)|
+ \frac{1}{M_1'} \left(1+ \sqrt{P_\delta}r \right)^2 
\Big| U'(r) + \frac{d-1}{r}U(r) \Big|
\\
& \qquad 
+ \frac{1}{M_2} \left(1+ \sqrt{P_\delta}r \right)^2 |\Theta (r)|
+ \frac{1}{M_2} (P_\delta r)^{-1}\left(1+ \sqrt{P_\delta}r \right)^2 |\Theta'(r)|.
\end{split}
\]
In the above, the constants are chosen such that 
$$
M_1 + M_1' + M_2 \ll 1 \qquad \mbox{and} \qquad M_1 < M_1';
$$
these constants will be determined more precisely shortly.

\begin{lem} 
\label{lemZvacuum}
Assume that $M_1' \log \frac{1}{P_\delta \delta^2}$ is sufficiently small.
Then $Z(\delta) < 1$ provided that
$$ 
\alpha \ll M_1, \qquad \alpha + \Theta_0 \ll M_2, \qquad \alpha^2 \ll M_2 P_\delta.
$$
Furthermore, if $Z(r) \leq 1$ for some $r>0$, then
$$
Z(r) \lesssim M_1 + \frac{P_\delta M_2}{M_1} + \frac{M_1}{M_1'} + \frac{\alpha}{M_1} + \frac{\Theta_0}{M_2} + \frac{M_1 M_1'}{P_\delta M_2}.
$$
\end{lem}

Finally, the asymptotic behavior of $U$, $\Theta$, and $P$ is established in the following lemma

\begin{lem} \label{asymptoticbehavvacuum}
There exists $P_\infty>0$, $U_\infty \in \mathbb{R}$,  $\Theta_\infty>0$ such that
\begin{align*}
P(r) = P_\infty + O \left( \frac{1}{r^2} \right), \; U(r) = \frac{U_\infty}{r} + O \left( \frac{1}{r^3} \right), \;\; \Theta(r) = \frac{\Theta_\infty}{r^2} + O\left( \frac{1}{r^4} \right).
\end{align*}
Furthermore, $\Theta(r)>0$ for any $r$ provided
$$
M_1 M_2 + \frac{M_1 M_1'}{P_\delta} \ll \Theta_0.
$$
Finally, $\Theta_\infty$ and $U_\infty$ can be expanded as
\begin{align*}
& U_{\infty} =  
\frac{2d(2\mu + \lambda)}{P_\delta} \alpha 
    \left[ 1 + O\Big( \alpha \log \frac{1}{P_\delta \delta^2} \Big) \right], \\
& \Theta_\infty = 
\frac{2(d-2)\kappa}{C_V P_\delta} \Theta_0 
\left[ 1+ O\left(\alpha \log \frac{1}{P_\delta \delta ^2}\right) 
 + O\left( \frac{\alpha^2 }{P_\delta \Theta_0} \right) 
 \right]. 
\end{align*}
\end{lem}

There remains to choose the constants $M_1$, $M_1'$, $M_2$, and to understand which range is allowed for parameters $\alpha, \delta, P_\delta, \Theta_0$, so that lemmas~\ref{lempsivacuum},~\ref{lemZvacuum} 
and~\ref{asymptoticbehavvacuum} 
apply. Let us summarize the requirements:
\begin{itemize}
\item We are interested in the perturbative regime where $\alpha,\delta,P_\delta,\Theta_0,M_1,M_1',M_2 \ll 1$. 
\item Local well posedness: requires $P_\delta \Theta_0 \ll \alpha$.
\item Control of $P$: requires $M_1' \log \frac{1}{P_\delta \delta^2} \ll 1$.
\item $Z(\delta) <1$: requires $\alpha \ll M_1$, $\alpha +  \Theta_0 \ll M_2$, $\alpha^2 \ll M_2 P_\delta$.
\item Bootstrap on $Z$: requires $M_1 + \frac{P_\delta M_2}{M_1} + \frac{M_1}{M_1'} + \frac{\alpha}{M_1} + \frac{\Theta_0}{M_2} + \frac{M_1 M_1'}{P_\delta M_2} \ll 1$.
\item Positivity of $\Theta$: requires $M_1 M_2 + \frac{M_1 M_1'}{P_\delta} \ll \Theta_0$.
\end{itemize}
First, we define $M_1$, $M_1'$, and $M_2$ as follows:
$$
M_1 = \Lambda \alpha, \qquad M_1' = \Lambda^2 \alpha, \qquad M_2 = \frac{\alpha}{P_\delta},
$$
where $\Lambda > 0$ is taken so big that $\frac{P_\delta M_2}{M_1} + \frac{M_1}{M_1'} + \frac{\alpha}{M_1} = \frac{3}{\Lambda} \ll 1$. Assuming that $\alpha,P_\delta, \Theta_0 \ll 1$, the remaining conditions reduce to the requirement that
$$
\alpha^2 \ll P_\delta \Theta_0 \ll \alpha, \qquad \mbox{and} \qquad \alpha \log \frac{1}{P_\delta \delta^2} \ll 1,
$$
which completes the proof.

\subsection{Local existence: proof of Lemma~\ref{lempsivacuum}}

\noindent \underline{Estimates on auxiliary functions.} 
Consider $(U,\Theta) \in B_{E^\delta}((\alpha r, \Theta_0),\frac{\alpha}{2})$, which implies in particular that $U(r) \geq \frac{\alpha}{2} r$. We derive various estimates on $V$, $W$, $Z$, $P$, $F_U$ and $F_\Theta$. First observe that, for $r < \delta$,
$$
V(\delta) - V(r) = \int_r^\delta \frac{U' + \frac{d-1}{r_1} U}{\frac{1}{2} r_1 - U} \, dr_1 \geq d\alpha \log \left( \frac{\delta}{r} \right).
$$
As a consequence,
$$
P(r) \leq P_\delta \left( \frac{r}{\delta} \right)^{d\alpha} .
$$
It is then easy to see that
\begin{align*}
|W| + |Z| \lesssim P_\delta \left( \frac{r}{\delta} \right)^{d\alpha} r^2 
\qquad \mbox{and} \qquad |W'| + |Z'| \lesssim P_\delta \left( \frac{r}{\delta} \right)^{d\alpha} r.
\end{align*}
Furthermore,
$$
|V'(r)| \lesssim \frac{\alpha}{r} \qquad \mbox{and} \qquad |P'(r)| \lesssim P_\delta \frac{\alpha}{r} \left( \frac{r}{\delta} \right)^{d\alpha}.
$$
We deduce from the above that
$$
|\widetilde{F}_U - d(2\mu + \lambda) \alpha| \lesssim 
P_\delta (\alpha + \Theta _0)
 \qquad \mbox{and} \qquad 
 |F_\Theta(r)| \lesssim 
 (P_\delta \alpha  \Theta _0 + \alpha^2 )r.
$$

\bigskip 

\noindent \underline{Estimates on differences of auxiliary functions.} Consider two elements $(U_1,\Theta_1)$ and $(U_2,\Theta_2)$ of 
$B_{E^\delta}((\alpha r, \Theta_0),\frac{\alpha}{2})$, with associated functions $V_1$, $V_2$, etc... We will denote $D$ their distance:
$$
D = \| (U_1,\Theta_1) - (U_2,\Theta_2) \|^\delta. 
$$
To start with, it is obvious that
$$
|\Theta_1(r) - \Theta_2(r)| \lesssim D r^2.
$$
Next observe that $|V_1' - V_2'| \lesssim \frac{D}{r}$, which implies $|[V_1(\delta) - V_1(r)] - [V_2(\delta) - V_2(r)]| \lesssim D \log \left( \frac{\delta}{r} \right)$. Using successively the inequalities $|e^x - e^y| \leq \max(e^x,e^y)|x-y|$ and $\sup_{0<t<1} t^\sigma \log \frac{1}{t} \sim \frac{1}{\sigma}$, this gives
$$
|P_1(r) - P_2(r)| \lesssim P_\delta \left( \frac{r}{\delta} \right)^{d\alpha} D \log\left( \frac{\delta}{r} \right) \lesssim  D \frac{P_\delta}{\alpha} . 
$$ 
This leads to 
\begin{align*}
& |W_1(r) - W_2(r)| + |Z_1(r) - Z_2(r)| \lesssim  D\frac{P_\delta }{\alpha} r^2, \\
& |\widetilde{F}_{U_1}(r) - \widetilde{F}_{U_2}(r)| \lesssim 
D \left[ P_\delta  + \frac{P_\delta\Theta_0}{\alpha} \right],
 \\
& |F_{\Theta_1}(r) - F_{\Theta_2}(r)| \lesssim 
D \left[ P_\delta \Theta_0 + \alpha \right] r.
\end{align*}

\bigskip 

\noindent \underline{Estimate on $\Psi(U_1,\Theta_1) - \Psi (U_2,\Theta_2)$.} Let $(\widetilde{U}_i,\widetilde{\Theta}_i) := \Psi(U_i,\Theta_i)$. Then
\begin{align*}
|\widetilde U_1(r) - \widetilde U_2(r)| \lesssim & r^{1-d} \int_0^r r_1^{d-1} |e^{W_1(r_1)  - W_1(r)} - e^{W_2(r_1)  - W_2(r)}| \widetilde F_{U_1}(r_1) \, dr_1 \\
& \qquad + r^{1-d} \int_0^r r_1^{d-1} e^{W_2(r_1)  - W_2(r_1)} |\widetilde F_{U_1}(r) - \widetilde F_{U_2}(r)| \,dr_1 \\
&
 \lesssim r^{1-d} \int_0^r r_1^{d-1} 
  \left[ P_\delta  + \frac{P_\delta\Theta_0}{\alpha} \right] D \,dr_1 
\\
& 
  \lesssim 
\left[ P_\delta + \frac{P_\delta\Theta_0}{\alpha} \right] D r.
\end{align*}
Similarly,
$$
|\widetilde U'_1(r) - \widetilde U'_2(r)| \lesssim \left[ P_\delta  + \frac{P_\delta\Theta_0}{\alpha} \right] D
$$
and finally,
\begin{align*}
& |\Theta_1(r) - \Theta_2(r)| \lesssim 
D \Big[\frac{P_\delta \Theta_0}{\alpha} + \alpha  \Big] r^2 , 
\\
& |\Theta_1'(r) - \Theta_2'(r)| \lesssim 
D \Big[\frac{P_\delta \Theta_0}{\alpha} + \alpha  \Big] r.
\end{align*}
As a conclusion,
$$
\| \Psi(U_1,\Theta_1) - \Psi (U_2,\Theta_2) \|^\delta \lesssim \left[ P_\delta + \alpha + \frac{P_\delta \Theta_0}{\alpha} \right] \| (U_1,\Theta_1) - (U_2,\Theta_2) \|^\delta.
$$
Therefore, $\Psi$ acts as a contraction provided $P_\delta + \alpha + \frac{P_\delta \Theta_0}{\alpha} \ll 1$. There remains to see that it stabilizes $B_{E^\delta}((\alpha r,\Theta_0),\alpha/2)$; this will be achieved in the following point.

\bigskip 

\noindent \underline{The first iterate.} Denoting $(\widehat{U},\widehat{\Theta}) = \Psi(\alpha r,\Theta_0)$, it follows that
\begin{align*}
& |\widehat{U} - \alpha r| \lesssim r P_\delta (\alpha + \Theta_0) ,
\\
& |\widehat{U}' - \alpha | \lesssim P_\delta (\alpha + \Theta_0) ,
\\
& |\widehat{\Theta} - \Theta_0| \lesssim (\alpha^2 + \Theta_0 P_\delta) r^2 , \\
& |\widehat{\Theta}' | \lesssim (\alpha^2 + \Theta_0 P_\delta) r.
\end{align*}
Therefore,
$$
\| \psi(\alpha r,\Theta_0) - (\alpha r,\Theta_0) \|^\delta \lesssim P_\delta \Theta_0 + \alpha^2 + \alpha P_\delta.
$$
This is $\ll \alpha$ provided $\alpha + P_\delta + \frac{P_\delta \Theta_0}{\alpha}$ is sufficiently small.

\bigskip 

\noindent
\underline{H\"older continuity of $P$}. 
We discuss around $r = 0$ only. Observe that 
\[
\begin{split}
& 
\Big|U'(r) + \frac{d-1}{r} U(r) \Big| 
= \Big| \frac{1}{2\mu + \lambda } \widetilde F_U (r) -W'(r) U(r) \Big|, 
\\
& 
\frac{1}{2\mu + \lambda } \widetilde F_U (r) 
= d \alpha + O(r^{d\alpha}), 
\quad 
W'(r) U(r) = O(r^2), 
\\
& 
\frac{1}{2}r - U(r) 
= \Big( \frac{1}{2} - \alpha \Big) r + \alpha r - U(r) 
= \Big( \frac{1}{2} - \alpha \Big) r + O(r^{1+d\alpha}). 
\end{split}
\]
These give that for $r \leq \delta$
\[
e^{V(r) - V(\delta) }
= \exp \Big( \int_\delta ^r 
 \frac{d\alpha}{\frac{1}{2}r_1 - \alpha r_1}  dr_1 
 + O(r^{1+d\alpha}) \Big) 
= \Big( \frac{r}{\delta}\Big) ^{\frac{2d\alpha}{1-2 \alpha}} 
\Big( 1 + O(r^{ 1+ d\alpha}) \Big),
\]
which proves $P(\cdot) \in C^{\frac{2d\alpha}{1 - 2\alpha}} (0,\delta)$.

\subsection{Global existence: proof of Lemma~\ref{lemZvacuum}}
\noindent
\underline{Estimate of $P$}. 
For $r \leq \delta$, it follows from the boundedness of $U$ that 
\[
V(\delta) - V(r) 
\geq \int_r^{\delta} 
\frac{\frac{\alpha}{2} + \frac{d-1}{s}\cdot \frac{\alpha}{2} s }
     {\frac{1}{2}s - \frac{\alpha}{2} s} ds 
\geq \frac{d\alpha}{1 - \alpha} \log\frac{\delta}{r} \geq d\alpha \log\frac{\delta}{r}
\]
and similarly, using $\alpha \ll 1$,
$$
V(\delta) - V(r) \leq 4 \alpha d \log \frac{\delta}{r}.
$$
Therefore,
\[
\left(\frac{r}{\delta}\right)^{4d\alpha} P_\delta
\leq 
P(r) 
\leq \left(\frac{r}{\delta}\right)^{d \alpha} P_\delta
\quad \text{ for any } r \in [0,\delta]. 
\]
As for $r \geq \delta$,
\begin{align*}
|V(r) - V(\delta)| = \left| \int_\delta^r \frac{U' + \frac{d-1}{s} U}{\frac{1}{2} s - U}\,ds \right| \lesssim \int_\delta^r \frac{M_1'}{s ( 1 + \sqrt{P_\delta} s)^2} \, ds \lesssim M_1' \left[ \log \frac{1}{P_\delta \delta^2} + 1\right] \lesssim 1
\end{align*}
under the assumption that $M_1' \log \frac{1}{P_\delta \delta^2}$ is $O(1)$.
This implies that, for $r \geq \delta$,
$$
P_\delta \lesssim P(r) \lesssim P_\delta.
$$
Finally, we record that, for any $r>0$,
\begin{align*}
& V'(r) \lesssim \frac{M_1'}{r ( 1 + \sqrt{P_\delta} r)^2} , \\
& W(r) + Z(r) \lesssim P_\delta r^2, \\
& W'(r) + Z'(r) \lesssim P_\delta r.
\end{align*}

\bigskip 

\noindent \underline{The starting point: $r=\delta$.} We start by checking that $Z(\delta) \ll 1$. Regarding the first three summands in the definition of $Z$, it is clear from the fixed point argument as soon as we choose
$$
\alpha \ll M_1 \qquad \mbox{and} \qquad \alpha + \Theta_0 \ll M_2.
$$
Turning to the fourth and last summand, we need to show that $|\Theta'(\delta)| \ll M_2 P_\delta r$. Observe that $\Theta$ has the same derivative as
$$
(d-2)r^{-d+2} \int_0^r r_1 ^{d-3} [ e^{-Z(r) + Z (r_1)} - 1] \, dr_1 \Theta_0 - \frac{U^2}{2C_V}
+ \frac{r^{-d+2}}{\kappa} \int_0^r 
  r_1^{d-2} e^{-Z(r) + Z(r_1)} F_\Theta (r_1)
 \, dr_1 ,
$$
which we write $I + II + III$. One can then check that, on $[0,\delta]$,
\begin{align*}
& |I'| \lesssim P_\delta \Theta_0 r, \\
& |II'| \lesssim \alpha^2 r, \\
& |III'| \lesssim (P_\delta \alpha \Theta_0 + \alpha^2) r.
\end{align*}
Therefore, $Z(\delta) \ll 1$ provided
$$
\alpha^2 \ll M_2 P_\delta.
$$

\bigskip

\noindent 
\underline{Estimate of $U$}. 
From a similar argument to \eqref{0913-1},
\[
e^{-CP_\delta (r^2 - r_1 ^2)}
\lesssim e^{-W(r)+W(r_1)} 
\lesssim e^{-C^{-1} P_\delta (r^2 - r_1 ^2)}
\quad \text{if } 0 \leq r_1 \leq r .
\]
Furthermore,
\begin{align*}
& \widetilde F_U(r) \lesssim 
 \Big( \frac{P_\delta r^2}{(1+\sqrt{P_\delta} r)^4} +1 \Big) M_1^2
 + \frac{P_\delta M_2}{(1+\sqrt{P_\delta}r)^2}  + \alpha 
\lesssim M_1^2 + P_\delta M_2 + \alpha,  \\
& \widetilde F_U'(r) \lesssim \frac{M_1' M_1 + P_\delta M_2}{r} .
\end{align*}
It follows by~\eqref{oriole} that
\[
|U(r)| 
\lesssim \frac{r}{(1+\sqrt{P_\delta}r)^2} ( M_1^2 + P_\delta M_2 + \alpha).
\]
As for the derivative, the required estimate for small $r $ is 
easy, since differentiating directly gives
\begin{align*}
\Big| U'(r) + \frac{d-1}{r}U(r) \Big| & = \left| \frac{1}{2\mu + \lambda} \widetilde F_U(r) - W'(r) U(r) \right| \\
&\lesssim M_1^2 + P_\delta M_2 + M_1 + \alpha\quad \text{for any } r \leq \frac{1}{\sqrt{P_\delta}} . 
\end{align*}
In order to investigate for large $r \geq 1/ \sqrt{P_\delta}$, first notice that 
\begin{equation}\notag 
\begin{split} 
\Big|
     \partial_{r} 
     \Big( \frac{r^{d-1} \widetilde F_U (r)}{W'(r)} \Big)
\Big| 
\lesssim
&  r^{d-3}  
\frac{\alpha + M_1 + P_\delta M_2}{P_\delta} 
\left( 1 + \Big(\frac{r}{\delta}\Big) ^{-4d\alpha} \right) , 
\end{split}
\end{equation}
where there appears a slight sigularity $r^{-4d\alpha }$ around $r = 0$
because of $P$ in the denominator, compared with \eqref{0913-2}, 
while it does not change the behavior. Using the inequality $r^{1-d} \int_0^r e^{-CP_\delta (r^2 - r_1^2)} r_1^{d-3}\,dr_1 \lesssim \frac{1}{P_\delta r^3}$ for $r \geq \frac{1}{\sqrt{P_\delta}}$, we get by~\eqref{doublestar}
\begin{align*}
\Big| U'(r) + \frac{d-1}{r}U(r) \Big| & = \left| W'(r) \frac{r^{1-d}}{{2\mu + \lambda}} \int_0^r e^{-W(r) + W(r_1)} \partial_{r_1} \left( \frac{r_1^{d-1} F_U(r_1)}{W'(r_1)} \right) \, dr_1 \right| \\
& \lesssim P_\delta r \frac{1}{P_\delta r^3} \left(\frac{\alpha + M_1 + P_\delta M_2}{P_\delta}\right) \\
& \lesssim \frac{1}{P_\delta r^2} \left( \alpha + M_1 + P_\delta M_2  \right)
\quad \text{for } r \geq \frac{1}{\sqrt{P_\delta}}.
\end{align*}

\bigskip 

\noindent
\underline{Estimate of $\Theta$}. 
Observe that for $r \leq 1 / \sqrt{P_\delta}$, 
\[
F_\Theta (r) 
\lesssim   (  M_1 M_1' + P_\delta M_1 M_2) r,
\]
which proves that 
\begin{align*}
|\Theta (r)| &\lesssim \Theta_0 + \frac{M_1 M_1' }{P_\delta} + M_1 M_2
     \quad \text{for } r \leq \frac{1}{\sqrt{P_\delta}}. 
\end{align*} 
Furthermore, for any $r>0$,
\[
F_\Theta(r) 
\lesssim \frac{1}{r} \left( M_1 M_2 + \frac{M_1 M_1'}{P_\delta} \right),
\]
which gives by~\eqref{oriole}
\[
|\Theta(r)| 
\lesssim \frac{1}{P_\delta r^2} \left( \Theta_0 + M_1 M_2 + \frac{M_1 M_1'}{P_\delta}\right) . 
\]
Finally, if $ r \leq \frac{1}{\sqrt{P_\delta}} $,
\[
|\Theta ' (r)| 
\lesssim P_\delta r \left( \Theta_0 + M_1 M_2 + \frac{M_1 M_1'}{P_\delta}\right)
\]
while if $r \geq \frac{1}{\sqrt{P_\delta}}$,
\[
|\Theta ' (r)| 
\lesssim \frac{1}{r} \left( \Theta_0 + M_1 M_2 + \frac{M_1 M_1'}{P_\delta}\right).
\]
This completes the proof of the desired estimate on $Z$.

\subsection{Asymptotic behavior: proof of Lemma~\ref{asymptoticbehavvacuum}}
The fact that $P$ converges, while $\Theta \sim \frac{\Theta_\infty}{r^2}$ and $U \sim \frac{U_\infty}{r}$, can be proved as in the smooth case and will not be developed here.

Turning to the positivity of $\Theta$, we observe first that, by the estimates above, there exists 
$$
\Theta(r) = (d-2) r^{-d+2} \int_0^r r_1^{d-3} e^{-Z(r) + Z(r_1)} \,dr_1 \Theta_0 + O \left( (M_1 M_2 + \frac{M_1 M_1'}{P_\delta} ) \frac{1}{(1 + \sqrt{P_\delta} r)^2} \right).
$$
Since, by~\eqref{oriole},
$$
(d-2) r^{-d+2} \int_0^r r_1^{d-3} e^{-Z(r) + Z(r_1)} \,dr_1 \Theta_0 \gtrsim \frac{\Theta_0}{(1 + \sqrt{P_\delta} r)^2},
$$
this means that $\Theta > 0$ provided
$$
M_1 M_2 + \frac{M_1 M_1'}{P_\delta} \ll \Theta_0.
$$

\noindent \underline{Expansion of $\Theta_\infty$ and $U_\infty$.} 
Observe that 
\[
P_\infty 
= P_\delta \left[ 1+ O\Big( V(\infty) - V(\delta)\Big) \right] 
= P_\delta \left[1 + O\left(\alpha \log \frac{1}{P_\delta \delta ^2}\right)\right ].
\]
Since $W(r) = \frac{P_\delta}{4(2\mu+\lambda)} r^2 
[ 1+ O\big( \alpha \log \frac{1}{P_\delta \delta ^2})\big]$, 
as $r \to \infty$, 
\[
\frac{r^{-d+1}}{2\mu + \lambda} 
\int_0^r r_1 ^{d-1} e^{-W(r) + W(r_1)} d (2\mu + \lambda) \alpha \, dr_1 
=
\frac{2d (2\mu + \lambda)}{P_\delta}\alpha 
\left[ 1 + O\left(\alpha \log \frac{1}{P_\delta \delta ^2}\right)\right ]
\frac{1}{r},
\]
and 
\[
\frac{r^{-d+1}}{2\mu + \lambda} 
\int_0^r r_1 ^{d-1} e^{-W(r) + W(r_1)} 
 \Big( \int_0^{r_1} \frac{d-1}{r_2} PU^2 \, dr_2 \Big) \, dr_1 
= O\left( \frac{\alpha^2}{P_\delta r}\right).
\]
These prove  
\[
U(r) = \frac{U_\infty}{r} + O\left( \frac{1}{r^3} \right) 
\quad \text{ with }
U_\infty = \frac{2d (2\mu + \lambda)}{P_\delta} \alpha 
\left[ 1 + O\left(\alpha \log \frac{1}{P_\delta \delta ^2}\right)\right ]. 
\]
Similarly, since 
$Z(r) = \frac{C_V P_\delta}{4\kappa}r^2 
 \big[ 1 + O(\alpha \log \frac{1}{P_\delta \delta ^2}) \big]$, 
 as $r \to \infty$, 
\[
(d-2) r^{-d+2} \int_0^r r_1^{d-3} e^{-Z(r) + Z(r_1)} \,dr_1 \Theta_0
= \frac{2(d-2)\kappa}{C_V P_\delta} \Theta_0 
\left[1 +  O\left(\alpha \log \frac{1}{P_\delta \delta ^2}\right)\right ]
\frac{1}{r^2},
\]
and 
\[
\frac{r^{-d+2}}{\kappa} \int_0^r r_1^{d-3} e^{-Z(r) + Z(r_1)} F_\Theta \,dr_1 
= O\left( \frac{\alpha \Theta_0 }{P_\delta r^2} + \frac{\alpha ^2}{P_\delta ^2 r^2} \right).
\]
These give 
\[
\Theta_\infty = 
\frac{2(d-2)\kappa}{C_V P_\delta} \Theta_0 
\left[ 1+ O\left(\alpha \log \frac{1}{P_\delta \delta ^2}\right) 
 + O\left( \frac{\alpha^2 }{P_\delta \Theta_0} \right) 
 \right]. 
\]


\begin{thebibliography}{00}

\bibitem{CP} M. Cannone, F. Planchon, Self-similar solutions for Navier-Stokes equations in ${\bf R}^3$, \textit{Comm. Partial Differential Equations} 21 (1996), no. 1-2, 179--193.

\bibitem{Chemin} J.-Y. Chemin, Th\'{e}or\`emes d'unicit\'{e} pour le syst\`eme de Navier-Stokes
   tridimensionnel, \textit{J. Anal. Math.} 77 (1999), 27--50.

\bibitem{CD-2015}
N. Chikami, R. Danchin, 
On the well-posedness of the full compressible Navier-Stokes system in critical Besov spaces, 
\textit{J. Differential Equations} 258 (2015), no. 10, 3435--3467.

\bibitem{D-2001}
R. Danchin, 
Global existence in critical spaces for flows of compressible viscous and heat-conductive gases, 
\textit{Arch. Ration. Mech. Anal.} 160 (2001), no. 1, 1--39. 

\bibitem{F-2004}
E. Feireisl, 
\textit{Dynamics of viscous compressible fluids.} 
Oxford Lecture Series in Mathematics and its Applications, 26. Oxford University Press, Oxford, 2004. 

\bibitem{FNP-2001}
E. Feireisl, A. Novotn\'y, H. Petzeltov\'a, 
On the existence of globally defined weak solutions to the Navier-Stokes equations, 
\textit{J. Math. Fluid Mech.} 3 (2001), no. 4, 358--392. 

\bibitem{GJ-2006} Z. Guo, S. Jiang, Self-similar solutions to the isothermal compressible
   Navier-Stokes equations, \textit{IMA J. Appl. Math.} 71 (2006), 658--669.

 \bibitem{HS-2001}
D. Hoff, J. Smoller, 
Non-formation of vacuum states for compressible Navier-Stokes equations, 
\textit{Comm. Math. Phys.} 216 (2001), no. 2, 255--276.

\bibitem{HL-2018}
X. Huang, J. Li, , 
Global classical and weak solutions to the three-dimensional full compressible 
Navier-Stokes system with vacuum and large oscillations, 
\textit{Arch. Ration. Mech. Anal.} 227 (2018), no. 3, 995--1059. 

\bibitem{JM-2015} 
J. Jang, N. Masmoudi, 
Well-posedness of compressible Euler equations in a physical vacuum, 
\textit{Comm. Pure Appl. Math.} 68 (2015), no. 1, 61--111.

\bibitem{JM-2015-2} 
J. Jang, N. Masmoudi, 
\textit{Vacuum in gas and fluid dynamics.} 
Nonlinear conservation laws and applications, 315--329, IMA Vol. Math. Appl., 153, Springer, New York, 2011. 

\bibitem{JS} H. Jia, V. Sverak, Local-in-space estimates near initial time for weak solutions of the Navier-Stokes equations and forward self-similar solutions, \textit{Invent. Math.} 196 (2014), no. 1, 233--265.

\bibitem{KT} H. Koch and D. Tataru, Well-posedness for the Navier-Stokes equations, \textit{Adv. Math.} 157 (2001), no. 1, 22--35.

\bibitem{LCX-2013} T. Li, P. Chen, J. Xie, Self-similar solutions of the compressible flow in one-space
   dimension, \textit{J. Appl. Math.} (2013), Art. ID 194704, 5.

\bibitem{Lions} P.-L. Lions, \textit{Mathematical topics in fluid mechanics. Vol. 1. Incompressible models.} Oxford Lecture Series in Mathematics and its Applications, 3. Oxford Science Publications. The Clarendon Press, Oxford University Press, New York, 1996. 

\bibitem{Lions-1998} P.-L. Lions, \textit{Mathematical topics in fluid mechanics. Vol. 2. Compressible models.} Oxford Lecture Series in Mathematics and its Applications, 10. Oxford Science Publications. The Clarendon Press, Oxford University Press, New York, 1998. 

\bibitem{MN-1980}
A. Matsumura, T. Nishida, 
The initial value problem for the equations of motion of viscous and heat-conductive gases, 
\textit{J. Math. Kyoto Univ.} 20 (1980), no. 1, 67--104. 


\bibitem{QSD-2008} Y. Qin, X. Su, S. Deng, Remarks on self-similar solutions to the compressible
   Navier-Stokes equations of a 1D viscous polytropic ideal gas, \textit{Appl. Math. Sci. (Ruse)} 2 (2008), 1493--1506.


\bibitem{Xin-1998}
Z. Xin, 
Blowup of smooth solutions to the compressible Navier-Stokes equation with compact density,  
\textit{Comm. Pure Appl. Math.} 51 (1998), no. 3, 229--240.

\end{thebibliography}
\end{document}